\documentclass[a4paper,12pt]{amsart}
\pdfoutput=1 

\usepackage[top=3cm,bottom=3cm,outer=3cm,inner=2cm,marginpar=2.45cm]{geometry}
\usepackage[pdftex,destlabel,final,colorlinks=true]{hyperref}
\usepackage[abbrev]{amsrefs}

\usepackage[french,ngerman,english]{babel}
\usepackage[utf8]{inputenc}
\usepackage[T1]{fontenc}

\usepackage[all,pdf]{xy}
\renewcommand{\objectstyle}{\displaystyle}

\usepackage{enumitem}
\usepackage{mathtools}  
\usepackage[usenames,dvipsnames]{xcolor}
\usepackage{calc} 

\babeltags{de = ngerman}
\babeltags{fr = french}

\usepackage[
show-mario,  
]{marionotations}

\usepackage{bbm}     

\usepackage[Smaller]{cancel} 

\usepackage{breakurl}  

\usepackage{subfiles}

\newcommand{\defaultDimension}{n}

\newcommand{\defaultAmbientSpace}{X}

\newcommand{\defaultlcIndex}{\sigma}

\newcommand{\defaultcohDegree}{q}

\newcommand{\defaultlclocus}{D}

\newcommand{\defaultvphi}{\vphi_F}

\newcommand{\defaultpsi}{\psi_D}

\newcommand{\defaultMetric}{\omega}

\NewDocumentCommand{\logKX}{
  t{M} 
  o    
}{K_X \otimes D \otimes \IfNoValueTF{#2}{F \IfBooleanT{#1}{\otimes M}}{#2}}

\NewDocumentCommand{\vphilist}{
  D||{\vphi}           
  t{F}                 
  t{M}                 
  d()                  
  D<>{\defaultMetric}  
}{\IfBooleanTF{#2}{\vphi_F}{#1} \IfBooleanT{#3}{+\vphi_M}, \IfNoValueF{#4}{(#4),} #5}

\NewDocumentCommand{\Ltwo}{ 
  D//{\bullet,\bullet}      
  D<>{\defaultAmbientSpace} 
  s                         
  m                         
}{L^{#1}_{(2)}\paren{\IfBooleanF{#3}{#2;} #4}}

\NewDocumentCommand{\Harm}{ 
  t{'}                      
  D//{\defaultcohDegree}    
  D<>{\defaultAmbientSpace} 
  g                         
  t{,}                      
  G{\defaultvphi}           
  e{_}                      
}{\mathcal{H}^{\IfBooleanF{#1}{\defaultDimension,}#2}\IfNoValueF{#4}{\paren{#3;#4}}_{#6 \IfNoValueF{#7}{,#7}}}

\NewDocumentCommand{\lcIndex}{ 
  m  
  m  
  m  
}{#1\IfNoValueF{#2}{+#2}\IfNoValueF{#3}{-#3}}

\NewDocumentCommand{\lcData}{ 
  G{\defaultvphi}  
  O{\defaultpsi}   
  e{.}             
}{\paren{#1; \IfNoValueF{#3}{#3 \cdot} #2}}

\NewDocumentCommand{\lcdata}{ 
  s                
  d<>              
  G{\defaultvphi}  
  O{\defaultpsi}   
  e{.,}            
}{\newcommand{\datalist}{\IfNoValueF{#2}{#2,}#3,#4\IfNoValueF{#5}{,#5}\IfNoValueF{#6}{,#6}}
\IfBooleanTF{#1}{\datalist}{\paren{\datalist}}}

\newcommand{\spHbase}{\mathbb{H}}
\NewDocumentCommand{\spH}{ 
  D//{\defaultcohDegree}  
  t{M}                    
  m                       
}{\spHbase^{#1}\paren{\IfBooleanT{#2}{M\otimes}#3}}

\DeclareMathOperator{\lc}{lc} 
\NewDocumentCommand{\lcc}{ 
  D||{\defaultlcIndex}       
  e{+-}                      
  D<>{\defaultAmbientSpace}  
  t{'}                       
  D(){\defaultlclocus}       
}{\lc_{#4}^{\lcIndex{#1}{#2}{#3}}\IfBooleanTF{#5}{\paren{#6}}{\lcData}}

\NewDocumentCommand{\lcS}{  
  s                       
  D(){\defaultlclocus}    
  D||{\defaultlcIndex}    
  e{+-}                   
  d<>                     
  O{p}                    
}{\mathtt{\IfBooleanT{#1}{\rs} #2}^{\lcIndex{#3}{#4}{#5}}_{\IfNoValueF{#6}{#6,}#7}}

\NewDocumentCommand{\PRes}{ 
  O{}      
  d()      
}{\mathcal R_{#1}\IfNoValueF{#2}{\paren{#2}}}

\newcommand{\defidlof}[1]{\mathcal{I}_{#1}}  
\NewDocumentCommand{\mtidlof}{   
  O{}      
  D<>{#1}  
  m        
}{\multidl_{#2}\paren{#3}} 

\NewDocumentCommand{\residlof}{  
  D||{\defaultlcIndex}   
  e{+-}                  
  d<>                    
  s                      
}{\sheaf R_{\IfNoValueTF{#4}{}{#4,} \lcIndex{#1}{#2}{#3}}\IfBooleanF{#5}{\lcData}}

\NewDocumentCommand{\Adjidlof}{
  D||{\defaultlcIndex}       
  D<>{\defaultAmbientSpace}  
  D(){\defaultlclocus}       
  m                          
}{\operatorname{\mathit{Adj}}^{#1}_{\paren{#2,#3}}\paren{#4}}

\NewDocumentCommand{\aidlof}{
  D||{\defaultlcIndex}   
  e{+-}                  
  d<>                    
  s                      
}{\sheaf{J}_{\!\IfNoValueTF{#4}{}{#4,} \lcIndex{#1}{#2}{#3}}\IfBooleanF{#5}{\lcData}}

\NewDocumentCommand{\faidlof}{
  D||{\defaultlcIndex}   
  e{+-}                  
  t{/}                   
  D||{\defaultlcIndex}   
  e{+-}                  
}{\fracAidlof{\lcIndex{#1}{#2}{#3}}{\lcIndex{#5}{#6}{#7}}}

\NewDocumentCommand{\fracAidlof}{
  m                  
  m                  
  d<>                
  s                  
  G{\defaultvphi}    
  O{\defaultpsi}     
  e{.}               
}{\frac{
    \aidlof|#1|<#3>*\IfBooleanF{#4}{\lcData{#5}[#6].{#7}}
  }{
    \aidlof|#2|<#3>*\IfBooleanF{#4}{\lcData{#5}[#6].{#7}}
  }}

\NewDocumentCommand{\lcV}{ 
  D||{\defaultlcIndex}    
  D//{\defaultvphi}       
  d()                     
  e{^}                    
  O{\defaultpsi}          
}{\:d\operatorname{lcv}^{#1\IfNoValueF{#4}{,\paren{#4}}}_{\IfNoValueF{#3}{#3,}#2}\left[#5\right]}

\NewDocumentCommand{\Ohvol}{ 
  D//{\defaultvphi} 
  d()               
  O{\defaultpsi}    
}{\dvol_{\IfNoValueF{#2}{#2,}#1}\left[#3\right]}

\newcommand{\dvol}{\:d\vol}

\NewDocumentCommand{\lcDataNormSubscript}{
  d<>                   
  s                     
  G{\defaultvphi}       
  O{\defaultpsi}        
  e{.}                  
  D||{\defaultlcIndex}  
  e{+-}                 
}{\IfNoValueF{#1}{#1,}
  \IfBooleanF{#2}{#3, \IfNoValueF{#5}{#5 \cdot} #4,}
  \lcIndex{#6}{#7}{#8}}

\newcommand{\RTFsym}{\mathfrak{F}} 
\NewDocumentCommand{\RTF}{ 
  s          
  G{\RTFsym} 
  o          
  >{\SplitArgument{1}{,}} d<> 
  d||        
  D(){\eps}  
  t{,}       
}{%
  \begingroup%
    \newif\ifsmasht%
    \IfBooleanTF{#1}{\smashttrue}{\smashtfalse}%
    \newif\ifboolup%
    \booluptrue%
    \IfNoValueT{#3}{\IfNoValueT{#4}{\IfNoValueT{#5}{\boolupfalse}}}%
    \newcommand{\supsrptstr}{\IfNoValueF{#3}{#3}\IfNoValueF{#4}{\inner#4}\IfNoValueF{#5}{\abs{#5}^2}}
    \newcommand{\RTFvar}{#6}
    #2\RTFprocess
}

\NewDocumentCommand{\RTFprocess}{
  o                     
  d<>                   
  t{,}                  
  G{\defaultvphi}       
  O{\defaultpsi}        
  e{.}                  
  D||{\defaultlcIndex}  
  e{+-}                 
}{\newcommand{\subsrptstr}{%
    \IfNoValueTF{#1}{
    \IfNoValueF{#2}{#2,}
    \IfBooleanT{#3}{#4,#5,\IfNoValueF{#6}{#6,}}
    \lcIndex{#7}{#8}{#9}}{#1}}%
  \newcommand{\srptstr}{\cramped{{}^{\supsrptstr}%
      \ifboolup _
      \fi{\ifboolup\displaystyle\fi\paren{\RTFvar}%
          \ifboolup {\scriptstyle \subsrptstr} \else _{\subsrptstr} \fi%
        }}}%
  \ifboolup%
    \ifsmasht%
      \smash[t]{
        \raisebox{\depthof{$\srptstr$} * \real{0.3}}{$\srptstr$}%
      }%
    \else%
      \raisebox{\depthof{$\srptstr$} * \real{0.3}}{$\srptstr$}%
    \fi%
  \else%
    \srptstr%
  \fi%
  \endgroup%
}

\NewDocumentCommand{\mtlog}{O{e} d() D||{\defaultpsi}}{\log\!#1^{\paren{#2}}\abs{#3}}
\NewDocumentCommand{\slog}{O{e} D||{\defaultpsi}}{\log\abs{#1 #2}}
\NewDocumentCommand{\dlog}{O{e} D||{\defaultpsi}}{\mtlog[#1](2)|#2|}

\NewDocumentCommand{\logpole}{ 
  D||{\defaultpsi}       
  D//{\defaultlcIndex}   
  E{.^}{{e}{1+\eps}}     
  s                      
}{\abs{#1}^{#2} \IfBooleanTF{#5}{\slog[#3]|#1|}{\paren{\slog[#3]|#1|}^{#4}}}

\DeclareFontFamily{OMX}{MnSymbolE}{}
\DeclareSymbolFont{MnLargeSymbols}{OMX}{MnSymbolE}{m}{n}
\SetSymbolFont{MnLargeSymbols}{bold}{OMX}{MnSymbolE}{b}{n}
\DeclareFontShape{OMX}{MnSymbolE}{m}{n}{
    <-6>  MnSymbolE5
   <6-7>  MnSymbolE6
   <7-8>  MnSymbolE7
   <8-9>  MnSymbolE8
   <9-10> MnSymbolE9
  <10-12> MnSymbolE10
  <12->   MnSymbolE12
}{}
\DeclareFontShape{OMX}{MnSymbolE}{b}{n}{
    <-6>  MnSymbolE-Bold5
   <6-7>  MnSymbolE-Bold6
   <7-8>  MnSymbolE-Bold7
   <8-9>  MnSymbolE-Bold8
   <9-10> MnSymbolE-Bold9
  <10-12> MnSymbolE-Bold10
  <12->   MnSymbolE-Bold12
}{}
\DeclareMathDelimiter{\llangle}{\mathopen}%
{MnLargeSymbols}{'164}{MnLargeSymbols}{'164}
\DeclareMathDelimiter{\rrangle}{\mathclose}%
{MnLargeSymbols}{'171}{MnLargeSymbols}{'171}

\newcommand{\iinner}[2]{\left\llangle#1,#2\right\rrangle}

\NewDocumentCommand{\idxup}{ 
  m                  
  O{\defaultMetric}  
  s                  
}{\paren{#1}^{\mathrlap{\!\IfBooleanTF{#3}{\smash[t]{#2}}{#2}}\makebox[\maxof{\widthof{$#2$}-\widthof{$\!\omega$}}{0pt}]{}}}

\newcommand{\dbadj}{\dbar^{\smash{\mathrlap{*}\;\:}}}

\NewDocumentCommand{\dep}{t{;} d<> O{\nu} m}{#4\IfBooleanTF{#1}{_}{^}{\IfNoValueF{#2}{#2\:}(#3)}}

\NewDocumentCommand{\sm}{s m}{{#2}\IfBooleanTF{#1}{_}{^}\text{sm}}

\NewDocumentCommand{\idx}{ 
  O{i} 
  m    
  o    
  t{.} 
  t{,} 
  o    
  m    
}{{#1}_{#2} \IfNoValueF{#3}{#3}
  \IfBooleanT{#4}{\dotsm} \IfBooleanT{#5}{,\dots,}
  \IfNoValueF{#6}{#6} {#1}_{#7}}

\NewDocumentCommand{\rs}{ 
  s  
  m  
}{\IfBooleanTF{#1}{\smash[t]{\widetilde{#2}}}{\widetilde{#2}}}

\DeclareMathOperator{\mlc}{mlc} 
\newcommand{\Diff}{\operatorname{Diff}^*} 

\newcommand{\sect}[1][s]{\mathtt{#1}} 

\NewDocumentCommand{\cbn}{  
  D//{\defaultlcIndex_V}
  D||{\defaultlcIndex}
}{\mathfrak{C}^{#1}_{#2}} 
\NewDocumentCommand{\Iset}{  
  D||{\defaultlcIndex}    
  e{+-}                   
  O{\defaultlclocus}      
  d()                     
}{I^{\lcIndex{#1}{#2}{#3}}_{#4}\IfNoValueF{#5}{\paren{#5}}} 

\ifcsname defineNoThmInMarionotations\endcsname
  \relax
\else 
  \newtheorem{prop}{Proposition}[section]
  
  \newtheorem{thm}[prop]{Theorem}

  \newtheorem{conjecture}[prop]{Conjecture}

  \theoremstyle{remark}
  \newtheorem{remark}[prop]{Remark}

  \theoremstyle{definition}
  \newtheorem{definition}[prop]{Definition}
  \newtheorem{example}[prop]{Example}

  \numberwithin{equation}{section}

\fi

\allowdisplaybreaks  

\begin{document}

\newcommand{\titlestr}{%
  An application of adjoint ideal sheaves \\
  to injectivity and extension theorems%
}

\newcommand{\shorttitlestr}{%
  An application of adjoint ideal sheaves to injectivity and extension
  theorems%
}

\newcommand{\MCname}{Tsz On Mario Chan}
\newcommand{\MCnameshort}{Mario Chan}
\newcommand{\MCemail}{mariochan@pusan.ac.kr}

\newcommand{\YJname}{Young-Jun Choi}
\newcommand{\YJnameshort}{Young-Jun Choi}
\newcommand{\YJemail}{youngjun.choi@pusan.ac.kr}

\newcommand{\PNUAddressstr}{%
 Department of Mathematics, Pusan National University, 2, Busandaehak-ro 63beon-gil, Geumjeong-gu, Busan 46241, Republic of Korea
}

\newcommand{\subjclassstr}[1][,]{%
  32J25 (primary)#1  
  32Q15#1   
  14B05 (secondary)
}

\newcommand{\keywordstr}[1][,]{%
  $L^2$ injectivity#1
  vanishing theorem#1
  extension theorem#1
  adjoint ideal sheaf#1
  multiplier ideal sheaf%
}

\newcommand{\dedicatorystr}{%
In honor of Professor Bo Benrdtsson on the occasion of his 70th birthday
}

\newcommand{\thankstr}{%
This work was supported by the National Research Foundation of
Korea (NRF) grant funded by the Korea government (MSIT)
(Nos.~2023R1A2C1007227 and 2021R1A4A1032418)%
}

\title
[\shorttitlestr]
{\titlestr}

\dedicatory{\dedicatorystr}
 
\author[\MCnameshort]{\MCname}
\email{\MCemail}

\author{\YJname}
\email{\YJemail}
\address{\PNUAddressstr}

\thanks{\thankstr}
 
\subjclass[2020]{\subjclassstr}

\keywords{\keywordstr}

\def\del{\partial}
\def\we{\wedge}
\def\ov{\overline}
\newcommand{\pd}[2]{\frac{\partial#1}{\partial#2}}

\date{\today} 

\maketitle


The version of the analytic adjoint ideal sheaf introduced in
\cite{Chan_adjoint-ideal-nas} is regarded by the authors as an
algebraic manifestation of the properties exhibited when the
residue function developed in \cite{Chan&Choi_ext-with-lcv-codim-1}
and \cite{Chan_on-L2-ext-with-lc-measures}, together with the
corresponding techniques of residue computation, is applied to (germs
of) holomorphic functions in multiplier ideal sheaves.
In \cite{Chan&Choi_injectivity-I} (as well as the upcoming joint work
with Shin-ichi Matsumura), it is shown that the adjoint ideal sheaves
and the residue computation fit in so well in the setup for solving
Fujino's conjecture, an injectivity theorem for log-canonical (lc)
pairs on compact K\"ahler manifolds.
The authors believe that the machinery can be useful for some other
problems involving lc pairs.

In this note, some properties of adjoint ideal sheaves and
some results of residue computations used in
\cite{Chan&Choi_injectivity-I} are reviewed.
Using only these results, a new ``qualitative'' extension result is
obtained in Section \ref{sec:extension}, which is essentially the
existence of holomorphic extensions suggested in
\cite{Chan_Residue-fct-proceedings}*{Conj.~2.2.3}, except that the result here
comes with no estimates.

\subsection{Notation and Conventions}\label{subsec:notation}
Unless stated otherwise, the following notations (which follow most of
the conventions in \cite{Chan&Choi_injectivity-I} and
\cite{Chan_adjoint-ideal-nas}) are used throughout this paper.

\begin{itemize}
\item $X$ is a compact K\"ahler manifold of dimension $n$. 

\item $D=\sum_{i \in \Iset||}D_{i}$ is a simple normal crossing (snc)
  reduced divisor on $X$, in which each $D_i$ is an irreducible
  component.
  Then, around every point in $D$, there is an
  open set $V \Subset X$ inside some coordinate chart with a
  holomorphic coordinate system $(z_1,\dots,z_n)$ such that
  \begin{equation} \label{eq:D-cap-V}
    D\cap V = \set{z_1 z_2 \dotsm z_{\sigma_V} = 0}
    \quad\text{ for some integer } \sigma_V \leq n \; .
  \end{equation}

\item $F$ is a holomorphic line bundle equipped with a (possibly
  singular) potential $\vphi_F$ (of the curvature of a metric
  $e^{-\vphi_F}$) which is \emph{quasi-plurisubharmonic (quasi-psh)}
  and has \emph{neat analytic singularities} such that $D
  \cup P_F$, where $P_F :=\vphi_F^{-1}(-\infty)$, is an snc divisor
  and $P_F$ contains no irreducible components of $D$ (hence no lc
  centres of $(X,D)$).
  In this case, around any point in $P_F$, there is
  an open set $V \Subset X$ on which \eqref{eq:D-cap-V} holds for some
  $\sigma_V \in \set{0,\dots, n}$ and 
  \begin{equation} \label{eq:vphi_F-on-V}
    \res{\vphi_F}_V \equiv \smashoperator[r]{\sum_{k =\sigma_V +1}^n} b_k \log\abs{z_k}^2
    \mod \smooth_X(V)
    \quad\text{   for some constants } b_k \geq 0 \; .
  \end{equation}

\item $X^\circ := X\setminus P_F$ and $\omega$ is the \emph{complete}
  K\"ahler form on $X^\circ$ of form given in
  \cite{Chan&Choi_injectivity-I}*{\S 2.2 item (4)}.

\item $\mtidlof{\vphi} := \mtidlof[X]{\vphi}$
  is the multiplier ideal sheaf of a potential $\varphi$, which is
  given at each $x \in X$ by 
  \begin{equation*}
    \mtidlof{\vphi}_x
    :=\setd{f \in \holo_{X,x}}{
      \begin{aligned}
        &f \text{ is defined on a coord.~neighbourhood } V_x \ni x \vphantom{f^{f^f}} \\
        &\text{and}\;\;\abs f^2 e^{-\vphi}\in L^1(V_x)
      \end{aligned}
    } \; .
  \end{equation*}

\item $\sect_i$ is a canonical section of the irreducible component
  $D_{i}$ for each $i \in \Iset||$, which gives a canonical section
  $\sect_D := \prod_{i\in \Iset||} \sect_i$ of $D$.
  
\item $\psi_D := \phi_D -\sm\vphi_D \leq -1$ is a global function on
  $X$ bounded above by $-1$, in which $\phi_D :=\log\abs{\sect_D}^2$
  is the potential on $D$ induced from the canonical section $\sect_D$
  and $\sm\vphi_D$ is a fixed smooth potential on $D$.

\item For each $\sigma \in \{0,1,\dots, n\}$, $\lcc' :=\bigcup_{p \in
    \Iset} \lcS$ is the union of \emph{$\sigma$-lc centres of
    $(X,D)$}, i.e.~the $\sigma$-codimensional irreducible components
  of any intersections of irreducible components of $D$ (under the
  assumption $(X,D)$ being log-smooth and lc), indexed by $\Iset$.
  Here $\lcc|0|' := X$ and $\Iset|0|$ being a singleton are set for
  convenience.



\item $\Diff_{p}D$ is the effective divisor on $\lcS$ defined by the 
  adjunction formula 
  \begin{equation*}
    K_{\lcS} \otimes \Diff_{p}D = \parres{K_X \otimes D}_{\lcS}
  \end{equation*}
  such that the restriction of $\sect_{(p)}:=
  \smashoperator{\prod\limits_{i \in \Iset \colon D_i
      \not\supset \lcS}} \sect_i$ to $\lcS$ is a canonical section of
  $\Diff_{p}D$.

\item $\psi_{(p)} :=\phi_{(p)} -\sm\vphi_{(p)} \leq -1$ is a global
  function on $X$ bounded above by $-1$ defined analogously to
  $\psi_D$, in which $\phi_{(p)} :=\log\abs{\sect_{(p)}}^2$ is the
  potential induced from $\sect_{(p)}$ and $\sm\vphi_{(p)}$ is a fixed
  smooth potential on the corresponding line bundle.
  
  %
  %
\end{itemize}
%
%

An open set $V \subset X$ is said to be \emph{admissible with respect
  to $(\vphi_F,\psi_D)$} if $V$ is biholomorphic to a polydisc centred
at the origin under a holomorphic coordinate system $(z_{1}, z_{2},
\cdots, z_{n})$ such that the residue computation in
\cite{Chan&Choi_injectivity-I}*{Thm.~2.6.1} is valid.
In particular, both \eqref{eq:D-cap-V} and \eqref{eq:vphi_F-on-V} are
satisfied on an admissible open set $V$.

When an admissible set $V$ is considered, an index $p \in \Iset$ such
that $\lcS \cap V \neq \emptyset$ is interpreted as a permutation
representing a choice of $\sigma$ elements from the set
$\set{1,2,\dots,\sigma_V}$ such that
\begin{equation*}
  \lcS \cap V = \set{z_{p(1)} = z_{p(2)} = \dotsm =
    z_{p(\sigma)} = 0}
  \quad\text{ and }\quad
  \res{\sect_{(p)}}_V = z_{p(\sigma+1)} \dotsm z_{p(\sigma_V)}
\end{equation*}
(cf.~the definition of the set $\cbn$ in
\cite{Chan_adjoint-ideal-nas}*{\S 3.1}).
Note also that the defining ideal sheaf
$\defidlof{\lcc|\sigma+1|'}$ of $\lcc|\sigma+1|'$ is
locally generated by $\res{\sect_{(p)}}_V = z_{p(\sigma+1)} \dotsm
z_{p(\sigma_V)}$ for all $p \in \Iset|\sigma|$.

\section{Properties of adjoint ideal sheaves}
\label{sec:properties-of-aidl}


In \cite{Chan_adjoint-ideal-nas}, the first author introduces the
following version of analytic adjoint ideal sheaves.

\begin{definition}[\cite{Chan_adjoint-ideal-nas}*{Def.~1.2.1}] \label{def:adjoint-ideal-sheaves}
  Given $\vphi_F$ and $\psi_D$ described as in Section
  \ref{subsec:notation}, the adjoint ideal sheaf $\aidlof{\vphi_F}
  :=\aidlof<X>{\vphi_F}$
  of index $\sigma$ is given at each $x \in X$ by
  \begin{equation*}
    \aidlof{\vphi_F}_x
    =\setd{f\in \holo_{X,x}}{\exists~\text{open set } V_x \ni x \: , \;
      \forall~\eps > 0 \: , \; \frac{\abs f^2
        e^{-\vphi_F-\psi_D}}{\logpole} \in L^1\paren{V_x} } \; . 
  \end{equation*}
\end{definition}
According to \cite{Chan_adjoint-ideal-nas}*{Thm.~1.2.3}, under the
assumption that $\vphi_F$ and $\vphi_F+\psi_D$ have only neat analytic
singularities with snc, one has
\begin{equation}\label{E:adjoint_ideal_sheaf_neat_analytic_sing}
  \aidlof{\vphi_F} =\mtidlof{\vphi_F} \cdot \defidlof{\lcc+1'}
\end{equation}
for all integers $\sigma \geq 0$, which fit into the chain of natural
inclusions 
\begin{equation*}
  \mtidlof{\vphi_F+\phi_D} =\aidlof|0|{\vphi_F}
  \subset \aidlof|1|{\vphi_F} \subset \dotsm \subset
  \aidlof|\sigma_{\mlc}|{\vphi_F} =\mtidlof{\vphi_F} \; ,
\end{equation*}
where $\sigma_{\mlc}$ is the codimension of the minimal lc centres
(mlc) of $(X,D)$.
Note that $\aidlof{\vphi_F} =\aidlof|\sigma_{\mlc}|{\vphi_F}$ for all
integers $\sigma \geq \sigma_{\mlc}$.

Moreover, since $\vphi_F^{-1}(-\infty)$ contains no lc centres of
$(X,D)$, the analytic adjoint ideal sheaves fit into the residue short
exact sequence
\begin{equation}\label{eq:short-ext-seq-of-ideals}
  \xymatrix@R-0.5cm@C+0.3cm{
    {0} \ar[r]
    & {\sheaf{J}_{\sigma-1}(\varphi_F;\psi_D)} \ar[r]
    & {\sheaf{J}_{\sigma}(\varphi_F;\psi_D)} \ar[r]^-{\Res^\sigma}
    & {\sheaf{R}_{\sigma}(\varphi_F;\psi_D)} \ar[r]
    & {0}
  } \; ,
\end{equation}
in which the \emph{residue sheaf $\residlof$ of index $\sigma$} is
given in this case as
\begin{equation*}
  \sheaf{R}_{\sigma}(\varphi_F;\psi_D)
  = \bigoplus_{p \in \Iset} \paren{\Diff_p D}^{-1}
  \otimes \mtidlof<\lcS>{\vphi_F} \; .
\end{equation*}
Note that the adjunction formula implies that
\begin{equation}\label{E:adjunction_residue_sheaf}
  K_X\otimes D\otimes F \otimes \sheaf{R}_{\sigma}(\varphi_F;\psi_D)
  =
  \bigoplus_{p \in\Iset} K_{\lcS} \otimes \res F_{\lcS} \otimes
  \mtidlof<\lcS>{\vphi_F} 
\end{equation}
(see \cite{Chan_adjoint-ideal-nas}*{\S 4.2} for the precise
construction of $\residlof{\vphi_F}$).
After tensoring with $\logKX$, the \emph{residue morphism
  $\Res^\sigma$} can be described in terms of the Poincar\'e residue
map $\PRes[\lcS]$ given in \cite{Kollar_Sing-of-MMP}*{\S 4.18} as
follows.
The Poincar\'e residue map $\PRes[\lcS]$ from $X$ to each $\lcS$ is
uniquely determined after an orientation on the conormal bundle of
$\lcS$ in $X$ is fixed.
For an admissible open set $V \subset X$, 
one has $\lcc' \cap V = \bigcup_{p \in \Iset}
\lcS<V>$ (where $\lcS<V> := \lcS \cap V$, which is connected by the
definition of the admissible open set, and possibly empty) and $\lcS<V>
=\set{z_{p(1)} =z_{p(2)} =\dotsm =z_{p(\sigma)}=0}$ when non-empty. 
Under such coordinate system, a section $f $ of  $\logKX \otimes
\aidlof$ on $V \subset X$ can be written as
\begin{equation*}
  f = \;\;\smashoperator{\sum_{p \in \Iset \colon \lcS<V>
      \neq\emptyset}} \;\; dz_{p(1)} \wedge \dotsm \wedge dz_{p(\sigma)}
  \wedge g_p \:\sect_{(p)} 
  =\;\;\smashoperator[l]{\sum_{p \in \Iset \colon \lcS<V>
      \neq\emptyset}}
  \frac{dz_{p(1)}}{z_{p(1)}} \wedge \dotsm
  \wedge \frac{dz_{p(\sigma)}}{z_{p(\sigma)}}
  \wedge g_p \:\sect_D \quad\text{ on } V \; , 
\end{equation*}
where $g_p$ is a local holomorphic $(n-\sigma, 0)$-form on $V$.
It therefore follows that 
\begin{equation*}
  \PRes[\lcS](\frac{f}{\sect_D})  =\res{g_p}_{\lcS} \in
  K_{\lcS} \otimes \res F_{\lcS} \otimes \mtidlof<\lcS>{\vphi_F}
  \quad\text{ on } \lcS<V> 
\end{equation*}
under the assumption that the orientation on the conormal bundle of
$\lcS$ in $X$ on $V$ is given by $(dz_{p(1)}, dz_{p(2)}, \dots,
dz_{p(\sigma)})$. 
The residue morphism $\Res^\sigma$ is then (according to
\cite{Chan_adjoint-ideal-nas}*{\S 4.2}) given by 
\begin{equation*}
  \renewcommand{\objectstyle}{\displaystyle}
  \xymatrix@C+0.5cm@R-0.5cm{
    {\logKX \otimes \aidlof} \ar[r]^-{\Res^\sigma}
    \ar@{}[d]|*[left]+{\in} 
    & {\logKX \otimes \residlof}
    \ar@{}[d]|*[left]+{\in}
    \\
    *+<0.8cm,0cm>{f} \ar@{|->}[r]
    & {\paren{\res{g_p}_{\lcS}}_{\mathrlap{p\in\Iset}}
      \mathrlap{\hphantom{p\in\Iset} .}} 
  }
\end{equation*}

In what follows, write
\begin{equation*}
  \aidlof* := \aidlof \quad\text{ and }\quad
  \residlof* := \residlof
\end{equation*}
for convenience when there is no risk of confusion.

Here is an explicit example which may help to illustrate the notion
described above.
\begin{example} 
Let $X=\Delta^n$ and $D=D_1+D_2$, where $D_1$ and $D_2$ are
hyperplanes defined by $D_1 :=\set{z_1=0}$ and $D_2 :=\set{z_2=0}$. 
Suppose that $\varphi_F$ is smooth, which implies that 
$\mtidlof{\vphi_F} =\sheaf{O}_X$. 
Then, it is obvious that $\mathrm{lc}_X^1(D)=D_1\cup D_2$,
$\mathrm{lc}_X^2(D)=D_1\cap D_2$, and $\sigma_{\mlc}=2$.
By \eqref{E:adjoint_ideal_sheaf_neat_analytic_sing}, the adjoint ideal sheaves are given by
\begin{equation*}
	\sheaf{J}_0
	=
	\mathcal{I}_{\mathrm{lc}_X^1(D)}
	=
	\mathcal{I}_{D_1\cup D_2}
	=
	\langle z_1z_2 \rangle
        \; , \quad
	\sheaf{J}_1
	=
	\mathcal{I}_{\mathrm{lc}_X^2(D)}
	=
	\mathcal{I}_{D_1\cap D_2}
	=
	\langle z_1,z_2 \rangle
        \;\;\;\text{and}\;\;\;
        \aidlof|2|* =\holo_X \; , 
\end{equation*}
where $\sheaf{J}_\sigma=\aidlof{\vphi_F}$ for $\sigma=0,1$.
Moreover, the residue sheaves are given by
\begin{equation*}
  \sheaf{R}_1 =\paren{\Diff_{D_1} D}^{-1} \oplus \paren{\Diff_{D_2} D}^{-1}
  =\res{D_2^{-1}}_{D_1} \oplus \res{D_1^{-1}}_{D_2}
  \quad\text{ and }\quad
  \residlof|2|* =\holo_{D_1 \cap D_2}  
\end{equation*}
(where $D_i^{-1}$ is the dual of $D_i$ as a line bundle).
For any local section $f\in \logKX \otimes \sheaf{J}_1$, it can be
written as
\begin{equation*}
f=\paren{z_2g_1 +z_1g_2} \:dz_1\wedge\cdots\wedge dz_n \; .
\end{equation*}
for some $g_1,g_2\in\sheaf{O}_X$.
The residue morphism $\Res^1$ is then given by
\begin{equation*}
	\Res^1(f)
	=
	\paren{
            	\parres{g_1 \:dz_2\wedge\cdots\wedge dz_n}_{D_1},
		\;
		-\parres{g_2 \:dz_1\wedge dz_3\wedge\cdots dz_n}_{D_2}
    	}\;.
\end{equation*}
\end{example}

Most properties of the adjoint ideal sheaves are derived from the
residue computations developed in
\cite{Chan&Choi_ext-with-lcv-codim-1}*{Prop.~2.2.1} and
\cite{Chan&Choi_injectivity-I}*{Thm.~2.6.1} (see also
\cite{Chan_on-L2-ext-with-lc-measures}*{Prop.~2.2.1}).
The computations are applicable even to smooth $(n,q)$-forms.
To state the result, let $\holoform_X$ be the holomorphic cotangent
bundle of $X$ and $\conj\holoform_X$ be its complex conjugate.
Write $\conj\holoform_X^q :=\bigwedge^q \conj\holoform_X$.
Extend each of the Poincar\'e residue maps $\PRes[\lcS]$ (hence
$\Res^\sigma$) so that it sends sections of $K_X \otimes
\conj\holoform_X^q$ (which are identified with $(n,q)$-forms in
$\smform/n,q/$) to $K_{\lcS} \otimes \res{\conj\holoform_X^q}_{\lcS}$.
Let $D +P_F =\sum_{k \in \Iset||[D+F]} E_k$ be the decomposition of the polar
set $D +P_F$ into irreducible components and let $\sect[e]_k$ be a
canonical section of $E_k$ for each $k \in \Iset||[D+F]$.
Set
\begin{align*}
  \smooth_{X\, *}
  &:=\paren{\smooth_{X}\left[
    \frac{1}{\abs{\sect[e]_k}} \colon
    k \in \Iset||[D+F] 
    \right]}_{\text{b}}
    \qquad\paren{\sect[e]_k \text{ treated as a local defining function of }
    E_k} \\
  &:=\set{\text{locally bounded elements in the $\smooth_X$-algebra generated
    by } \frac{1}{\abs{\sect[e]_k}} \text{ for } k\in\Iset||[D+F]}
    \; .\footnotemark
\end{align*}%
\footnotetext{
  On an admissible open set $V$ under the holomorphic coordinate
  system $(z_1,\dots, z_n)$ such that $D\cap V =\set{z_1 z_2 \dotsm
    z_{\sigma_V} =0}$ and $P_F \cap V =\set{z_{\sigma_V+1} \dotsm
    z_{\tau_V} = 0}$, one has
  \begin{equation*}
    \smooth_{X \,*}(V)
    =\smooth_X(V)\left[e^{\pm \cplxi \theta_1}, \dots, e^{\pm \cplxi
        \theta_{\tau_V}} \right]
  \end{equation*}
  where $(r_k,\theta_k)$ is the polar coordinate system of the
  $z_k$-plane for $k=1,\dots,\tau_V$ in $V$, which is (almost) the
  same as the ad hoc definition of $\smooth_{X\, *}(V)$ given in 
  \cite{Chan&Choi_injectivity-I}*{\S 2.6} (in which
  $e^{\pm\cplxi\theta_{k}}$ for $k \geq \tau_V +1$ are also included
  in the set of generators of the algebra).
  The definition given here is independent of coordinates and its
  sheaf structure can be seen easily.
}%
Set also
\begin{equation*}
  \smform/p,q/[X \,*] :=\set{\text{$(p,q)$-forms with coefficients in }
    \smooth_{X\,*}}.
\end{equation*}
The following proposition recalls the main result of the
residue computations.

\begin{prop}[\cite{Chan&Choi_injectivity-I}*{Thm.~2.6.1}]
  \label{prop:residue-product-X-to-lcS}
  Let
  $f \in \smform/n,q/[X\,*]\paren{X; D\otimes F \otimes \aidlof*}$,
  one has
  \begin{align*}
    \lim_{\eps \tendsto 0^+} \eps \int_X
    \frac{\abs{f}_\omega^2 \:e^{-\phi_D-\vphi_F}}{\abs{\psi_D}^{\sigma
    +\eps}}
    &=\sum_{p \in \Iset} \frac{\pi^\sigma}{(\sigma-1)!}
      \int_{\lcS} \abs{\PRes[\lcS](\frac{f}{\sect_D})}_\omega^2 e^{-\vphi_F}
      <+\infty \; .
  \end{align*}
\end{prop}

\section{Injectivity theorem for lc pairs}
\label{sec:injectivity}

Koll\'ar's injectivity theorem (\cite{Kollar_injectivity}*{Thm.~2.2})
can be considered as a generalization of the Kodaira vanishing theorem
(see, for example, \cite{Esnault&Viehweg_book}*{Cor.~5.2} or
\cite{Lazarsfeld_book-I}*{Remark 4.3.8}).
It was first proved in \cite{Kollar_injectivity} via algebraic method
and later its generalization to the setting on compact K\"ahler
manifolds was proved via harmonic theory by Enoki (\cite{Enoki}).
The theorem is further generalized to lc pairs in the algebraic
setting using the theory of (mixed) Hodge
structure (see, for example, \cite{Esnault&Viehweg_book}*{\S 5},
\cite{Fujino_log-MMP}*{\S 5} and \cite{Ambro_injectivity}).
Fujino formulated the following conjecture, expecting that it could be
proved following the line of arguments in Enoki's proof.

\begin{conjecture}[Fujino's conjecture, \cite{Fujino_survey}*{Conj.~2.21}] \label{conj:Fujino-conj}
  Given $F$ as described in Section \ref{subsec:notation}, which is
  semi-positive in particular, suppose that there exists a holomorphic
  section $s$ of $F^{\otimes m}$ on $X$ for some positive integer $m$
  such that $s$ does not vanish identically on any lc centres of
  $(X,D)$.
  Then, the multiplication map induced by $\otimes s$,
  \begin{equation*}
    \cohgp q[X]{K_X \otimes D\otimes F} \xrightarrow{\;\otimes s\;}
    \cohgp q[X]{K_X \otimes D\otimes F^{\otimes (m+1)}} \; ,
  \end{equation*}
  is injective for every $q \geq 0$.
\end{conjecture}
In the spirit of Enoki's arguments, Fujino
(\cite{Fujino_injectivity}), Matsumura (\cite{Matsumura_injectivity})
and Gongyo--Matsumura (\cite{Gongyo&Matsumura}) prove the
corresponding injectivity theorem for Kawamata log-terminal (klt)
pairs while Matsumura also obtains in \cite{Matsumura_injectivity-lc}
the result for purely log-terminal (plt) pairs (see these references
for more information on the development of the injectivity theorem).

Inspired by the idea of Matsumura for the plt case in
\cite{Matsumura_injectivity-lc}, the authors introduce the use of
adjoint ideal sheaves to the study of Fujino's conjecture in
\cite{Chan&Choi_injectivity-I}, which leads to the scheme of a proof
by induction on the codimension $\sigma$ of the lc centres of $(X,D)$
(see \eqref{eq:commut-diagram_sing-Fujino-conj} and Theorem
\ref{thm:induction-on-Fujino-conj}).
Moreover, a slightly more general statement, namely, Fujino's
conjecture with multiplier ideal sheaves, will follow naturally under
this framework.
It will be shown in our upcoming article with Shin-ichi Matsumura
that a proof of Fujino's conjecture in the spirit of Enoki's arguments
can be obtained by following this scheme.

A solution to Fujino's conjecture is announced recently by Cao and P\u
aun (\cite{Cao&Paun_LC-inj}), in which a different approach is used. 
Their approach, instead of working within the framework of $L^2$ theory,
makes use of a version of the Hodge decomposition on currents induced
via duality from the Hodge decomposition of $L^2$ forms with respect
to a K\"ahler form with conic singularities (along $D$).

The rest of this section is devoted to reviewing the statements in
\cite{Chan&Choi_injectivity-I} which bring adjoint ideal sheaves and
certain residue computations into play in the study of Fujino's conjecture.



For the sake of convenience, set
\begin{equation*}
  \spH{\sheaf F} :=\cohgp q[X]{K_X\otimes D\otimes
                                   F\otimes \sheaf F}
\end{equation*}
for any sheaf $\sheaf F$ on $X$ and any integer $q =0, \dots, n$.
Let $s \in \cohgp 0[X]{M}$ be a non-trivial holomorphic section of a
line bundle $M$ (where $M := F^{\otimes m}$ in Fujino's conjecture)
and $\vphi_M$ a potential on $M$ with the same regularity property as
$\vphi_F$ described in Section \ref{subsec:notation} such that $\sup_X
\abs s_{\vphi_M}^2 < \infty$.
From the residue short exact sequence
\eqref{eq:short-ext-seq-of-ideals} and the multiplication map
\begin{equation*}
  K_X\otimes D\otimes F\otimes \aidlof{\vphi_F}
  \xrightarrow{\otimes s\;}
  K_X\otimes D\otimes F\otimes M\otimes \aidlof{\vphi_F +\vphi_M} \; ,
\end{equation*}
one obtains the following commutative diagram of cohomology groups:

\begin{equation} \label{eq:commut-diagram_sing-Fujino-conj}
  \begin{gathered}
    \xymatrix@R=0.85cm@C+0.75cm{
      {\vdots} \ar[d]
      & {\vdots} \ar[d]
      & {\vdots} \ar[d]
      \\
      {\spH{\aidlof-1*}} \ar[d] \ar@{=}[r]
      & {\spH{\aidlof-1*}} \ar[r]^-{\otimes s}
      \ar[d]_-{\iota_{\sigma-1}} \ar[dr]|-*+{\mu_{\sigma-1}}
      & {\spH M{\aidlof-1{\vphi_F+\vphi_M}}} \ar[d]
      \\
      {\spH{\aidlof*}} \ar[d] \ar[r]^-{\iota_{\sigma}}
      \ar@/_1.9pc/[rr]|(.65)*+{\mu_{\sigma}}
      & {\spH{\aidlof|\sigma_{\mlc}|*}}
      \ar[d]|(.42)*+<3pt>{ }
      \ar[r]^-{\otimes s}
      & {\spH M{\mtidlof{\vphi_F+\vphi_M}}}
      \ar[d]
      \\
      {\spH{\residlof*}} \ar[d] \ar[r]^-{\tau_\sigma}
      \ar@/_1.65pc/[rr]+<-39pt,-15pt>|(.67)*+{\nu_\sigma}
      & {\spH{\frac{\aidlof|\sigma_{\mlc}|*}{\aidlof-1*}}} \ar[d]|(.5)*+<3pt>{}
      \ar[r]^-{\otimes s}
      & {\spH M{\frac{\mtidlof{\vphi_F+\vphi_M}}{\aidlof-1{\vphi_F+\vphi_M}}}} \ar[d] \\
      {\vdots} & {\vdots} & {\vdots} }
  \end{gathered}
\end{equation}
Note that the columns are all exact.
The 
map $\iota_\sigma$ is induced from the
natural inclusion $\aidlof* \subset \aidlof|\sigma_{\mlc}|* =\mtidlof{\vphi_F}$, while
the horizontal maps on the right-hand-side are induced from the
multiplication map $\otimes s$.
Each homomorphism of $\mu_\sigma$'s and $\nu_\sigma$'s is the
composition of the maps on the corresponding row.

Through a simple diagram-chasing, one sees that, for each $\sigma \geq
1$, if the homomorphisms $\mu_{\sigma-1}$ and $\nu_\sigma$ satisfy
$\ker\mu_{\sigma-1} =\ker\iota_{\sigma-1}$ and $\ker\nu_\sigma
=\ker\tau_\sigma$ respectively, then it follows that $\ker\mu_{\sigma}
=\ker\iota_\sigma$.
One then obtains the following theorem via induction.
\begin{thm} \label{thm:induction-on-Fujino-conj}
  If one has $\ker\mu_0 =\ker\iota_0$ and $\ker\nu_{\sigma}
  =\ker\tau_{\sigma}$ for $\sigma =1, \dots, \sigma_{\mlc}$, then
  $\mu_{\sigma_{\mlc}}$ is injective (as $\iota_{\sigma_{\mlc}}$ is the
  identity map).
  In particular, when $M :=F^{\otimes m}$ and when $\vphi_M :=
  m\vphi_F$ is smooth, Fujino's conjecture holds true under these
  given assumptions.
\end{thm}

In \cite{Chan&Choi_injectivity-I}, the assumption $\ker\mu_0
=\ker\iota_0$ is proved when $\varphi_F$ and $\varphi_M$ have neat
analytic singularities with a suitable positivity assumption.
More precisely, 

\begin{thm}[\cite{Chan&Choi_injectivity-I}*{Thm.~1.2.1}] 
  \label{thm:main-result}
  Suppose $(X,\phi_D,\vphi_F)$ is as described in Section \ref{subsec:notation}
  and let $\vphi_M$ be a quasi-psh potential on $M$ having the same regularity
  properties as $\vphi_F$ described in Section \ref{subsec:notation}
  such that
  \begin{itemize}
  \item $\ibddbar\vphi_F \geq 0$ and $\ibddbar\vphi_M
    \leq C\ibddbar\vphi_F$ on
    $X$ in the sense of currents for some constant $C > 0$ (where
    $\ibar := \ibardefn
    \;$\ibarfootnote), and 
  \item 
    $P_F\cup P_M \cup D$ has only snc (where $P_M :=\vphi_M^{-1}(-\infty)$).
  \end{itemize}
  Suppose further that there exists a non-trivial holomorphic section $s \in
  \cohgp 0[X]{M}$ such that $\sup_X \abs{s}_{\vphi_M}^2 < \infty$.
  Then, given the commutative diagram
  \begin{equation*}
    \xymatrix@R=0.6cm{
      {\cohgp q[X]{K_X\otimes D\otimes F\otimes \mtidlof{\phi_D+\vphi_F}}}
      \ar[d]_-{\iota_0} \ar@{}@<8ex>[d]|*+{\circlearrowleft} \ar[dr]^-{\mu_0}
      \\
      {\cohgp q[X]{K_X\otimes D\otimes F\otimes \mtidlof{\vphi_F}}}
      \ar[r]^-{\otimes s}
      &
      {\cohgp q[X]{K_X\otimes D\otimes F\otimes M\otimes
        \mtidlof{\vphi_F+\vphi_M}} \; ,}
    }
  \end{equation*}
  in which $\iota_0$ is induced from the inclusion
  $\mtidlof{\phi_D+\vphi_F} \subset \mtidlof{\vphi_F}$,
  one has $\ker\mu_0 =\ker\iota_0$ for every $q \geq 0$.
\end{thm}

This theorem with smooth $\varphi_F$ and $\varphi_M$ is first proved
by Matsumura by smoothing the canonical potential $\phi_D$ of $D$
(\cite{Matsumura_injectivity-lc}*{Thm.~3.9}) when applying the
Bochner--Kodaira formula with a (smooth) K\"ahler form on $X$.
Instead of smoothing out the lc singularities of $\phi_D$,
\cite{Chan&Choi_injectivity-I} considers the (twisted)
Bochner--Kodaira formula (hence, in particular, the formal adjoint
operator $\dfadj$ of $\dbar$) with respect to the \emph{complete}
K\"ahler form $\omega$ on $X^\circ = X \setminus \paren{P_F \cup P_M}$
(or simply the K\"ahler form on $X$ when $\vphi_F$ and
$\vphi_M$ are smooth) and the potential $\phi_D+\varphi_F
\:\paren{+\varphi_M}$ (which is singular along general points of $D$
on $X^\circ$).
Such Bochner--Kodaira formula is shown to be applicable to the harmonic
$L^2$ forms with respect to the data $(\phi_D,\vphi_F,\vphi_M,\omega)$
on $X^\circ$ by making use of a family of smooth cut-off functions
$\seq{\theta_\eps}_{\eps > 0}$ supported away from $D$ in $X^\circ$
such that $\theta_\eps \ascendsto 1$ almost everywhere in $X^\circ$ as
$\eps \descendsto 0$.
The proof (see \cite{Chan&Choi_injectivity-I}*{Prop.~3.2.5 and
  Cor.~3.2.6}) takes advantage of the residue computation in
Proposition \ref{prop:residue-product-X-to-lcS}.
A flavour of the computation can also be found in the following idea of
proof of Theorem \ref{thm:main-result} (in the evaluation of
$\iinner{su}{\dbar v_{(\infty)}}_{\vphi+\vphi_M, \omega}$).

The outline of the proof of Theorem \ref{thm:main-result} is given
without the treatment to the singularities on $\varphi_F$ and
$\varphi_M$, as the proof for the case when they are smooth already
illustrates the use of the residue computation.
For the full treatment for singular $\vphi_F$ and $\vphi_M$, see
\cite{Chan&Choi_injectivity-I}*{\S 3.3}.

\begin{proof}[Idea of the proof of Theorem \ref{thm:main-result}]
  Note that
  any class in $\cohgp q[X]{\logKX \otimes \mtidlof{\vphi}}$ is
  represented by a unique element in the space of harmonic $L^2$ forms
  $\mathcal{H}^{n,q}_{\varphi,\omega}$ on $X^\circ$ with respect to
  $\varphi=\phi_D+\vphi_F$ and $\omega$ (see
  \cite{Matsumura_injectivity}*{Prop.~5.5} and
  \cite{Matsumura_injectivity-lc}*{Prop.~2.8} for a proof).
  It is enough to show that, for any 
  \begin{equation*}
    [u]\in\cohgp q[X]{\logKX \otimes \mtidlof{\varphi}} 
    \quad\text{ with }\quad
    u\in\mathcal{H}^{n,q}_{\varphi,\omega} \; ,
  \end{equation*}
  if $[u]\in(\ker\iota_0)^\perp\cap\ker\mu_0$ (in which the orthogonal
  complement $\paren\cdot^{\perp}$ is taken with respect to the $L^2$
  norm $\norm\cdot_{\vphi,\omega}$ on $\Harm,{\vphi}_\omega$), then $u=0$.

  Since $0=[su]\in\cohgp q[X]{\logKX M \otimes
    \mtidlof{\vphi_F+\vphi_M}}$, $su =\dbar v$ for some $v\in
  L^{n,q-1}_{(2)}(D\otimes F\otimes
  M)_{\sm\varphi_D+\vphi_F+\vphi_M,\omega}$ which can be decomposed
  into $v = v_{(2)} + v_{(\infty)}$ such that
  \begin{align*}
    v_{(2)} &\in L^{n,q-1}_{(2)}\paren{D\otimes F\otimes
              M}_{\phi_D+\vphi_F+\vphi_M,\omega}  \quad\text{ and} \\
    v_{(\infty)} &\in \smform/n,q-1/\paren{X; D\otimes F\otimes M
                   \otimes \mtidlof{\vphi_F+\vphi_M}}
  \end{align*}
  via the $L^2$ Dolbeault isomorphism (see
  \cite{Chan&Choi_injectivity-I}*{Lemma 3.2.1}).
  Furthermore, the positivity $\ibddbar\vphi_F \geq 0$, together with
  the refined hard Lefschetz theorem (see
  \cite{Matsumura_injectivity-lc}*{Thm.~3.3} or
  \cite{Chan&Choi_injectivity-I}*{Thm.~2.5.1}), implies that $*_\omega
  u$ is holomorphic on $X$, which then implies that $\frac u{\sect_D}$
  is smooth on $X^\circ$.
  From
  \begin{equation*}
    \norm{su}_{\vphilist M}^2
    =\iinner{su}{\dbar v}_{\vphilist M}
    =\iinner{su}{\dbar v_{(2)}}_{\vphilist M}
    +\iinner{su}{\dbar v_{(\infty)}}_{\vphilist M} \; ,
  \end{equation*}
  it suffices to show that both of the inner products vanish, which
  then results in $su = 0$, hence $u=0$.

  \begin{description}[
    leftmargin=0pt, itemindent=*, align=left, itemsep=1.5ex, topsep=2ex]

  \item[The 1st term]
    $\iinner{su}{\dbar v_{(2)}}_{\vphilist M}
    =\iinner{\dfadj_{\vphi_M} \paren{su}}{v_{(2)}} =0$,
    where $\dfadj_{\vphi_M}$ is the formal adjoint operator of $\dbar$ 
    with respect to $\varphi+\varphi_M$ and $\omega$.
    \\
        
    From the fact that $\frac{u}{\sect_D}$ is smooth on $X^\circ$ and
    by making use of the family of cut-off functions
    $\seq{\theta_\eps}_{\eps > 0}$ supported away from $D$, the
    Bochner--Kodaira formula with respect to the potential $\vphi$ and
    metric $\omega$ is valid for $u$ (see
    \cite{Chan&Choi_injectivity-I}*{Prop.~3.2.5}). 
    Together with the positivity $\ibddbar\vphi_F \geq 0$ on $X$, one obtains
    \begin{equation} \label{eq:nable-u=0_ibddbar-vphi_F-u=0}
      \nabla^{(0,1)}u = 0 \quad\text{ and }\quad
      \idxup{\ibddbar\vphi_F}\ptinner{u}{u}_{\omega} = 0 
      \quad\text{ on } X^\circ
      \qquad\text{(\cite{Chan&Choi_injectivity-I}*{Prop.~3.2.5}),}
    \end{equation}
    where $\nabla^{(0,1)}$ is the covariant differential operator of type $(0,1)$ 
    induced from the Chern connection.
    Consequently, these results, together with the Bochner--Kodaira
    formula with respect to $\vphi+\vphi_M$ and $\omega$, lead to
    \begin{equation*} 
      \dfadj_{\vphi_M} \paren{su} = 0 \quad\text{ and }\quad
      su \in \Dom \dbadj_{\vphi_M}
      \qquad\text{(\cite{Chan&Choi_injectivity-I}*{Cor.~3.2.6}),}
    \end{equation*}
    where
    $\dfadj_{\vphi_M}$ is the formal adjoint of $\dbar$ with respect
    to the $L^2$ norm $\norm\cdot_{\vphi+\vphi_M, \omega}$ and
    $\dbadj_{\vphi_M}$ is the corresponding Hilbert space adjoints on
    the $L^2$ space $L_{(2)}^{n,\bullet}(D\otimes F\otimes
    M)_{\varphi+\varphi_M,\omega}$.
    Since $v_{(2)}$ is $L^2$ with respect to
    $\norm\cdot_{\vphi+\vphi_M, \omega}$ and lies inside $\Dom\dbar$,
    one can apply the integration by parts and the inner product
    under consideration vanishes as desired. 


  \item[The 2nd term] $\iinner{su}{\dbar v_{(\infty)}}_{\vphilist M}=0$.
    \\
    
    Note that $v_{(\infty)}$ is $L^2$ with respect to the potential
    $\vphi_F+\vphi_M$ but not necessarily to $\phi_D+\vphi_F+\vphi_M$.
    To overcome this difficulty, the cut-off function $\theta_\eps$
    (as described before) is introduced (and defined explicitly here)
    to facilitate the integration by parts.
    Let $\theta \colon [0,\infty) \to [0,1]$ be a smooth non-decreasing
    cut-off function such that $\res\theta_{[0,\frac 12]} \equiv 0$ and
    $\res\theta_{[1,\infty)} \equiv 1$.
    For $\eps \geq 0$, set $\theta_\eps := \theta \circ \frac
    1{\abs{\psi_D}^\eps}$ and $\theta'_\eps := \theta' \circ \frac
    1{\abs{\psi_D}^\eps}$ for convenience (where $\theta'$ is the
    derivative of $\theta$).
    Note that both $\theta_{\eps}$ and $\theta'_{\eps}$ have compact
    supports inside $X \setminus D$ for $\eps > 0$.
    One also has $\theta_\eps \nearrow 1$ almost everywhere in $X$ as $\eps
    \searrow 0$. As $\eps \tendsto 0^+$, one has
    \begin{align*}
      \iinner{su}{\dbar v_{(\infty)}}_{\vphilist M}
      &\leftarrow \iinner{su}{ \theta_\eps \dbar v_{(\infty)}}_{\vphilist M} \\
      &=
        \iinner{su}{ \dbar \paren{\theta_\eps v_{(\infty)}}}_{\vphilist M}
        -\iinner{su}{ \frac{\eps \theta'_\eps}{\abs{\psi_D}^{1+\eps}}
        \dbar\psi_D \wedge v_{(\infty)}}_{\vphilist M} 
      \\
      &=
        \cancelto{0}{\iinner{\dfadj_{\vphi_M} \paren{su}}{ \theta_\eps  v_{(\infty)}}}_{\vphilist M}
        -
        \eps \:\iinner{\idxup{\diff\psi_D}[\omega] \ctrt su}
        {\:\frac{\theta'_\eps v_{(\infty)}}{\abs{\psi_D}^{1+\eps}} 
        }_{\vphilist M} \; .
    \end{align*}

    Hinted by the argument $\idxup{\diff\psi_D} \ctrt u$ in the inner
    product, one is led to consider the twisted Bochner--Kodaira
    formula with respect to $\vphi$, $\omega$ and the twisting factor
    $\abs{\psi_D}^{1-\eps}$ for $\eps > 0$, given by
    \begin{align*}
      &~
        \int_{X^\circ} \abs{\dbar\zeta}_{\vphilist}^2 \:\abs{\psi_D}^{1-\eps}
        +\int_{X^\circ} \abs{\dfadj\zeta}_{\vphilist}^2 \:\abs{\psi_D}^{1-\eps}
      \\ 
      =&\!\!\!
         \begin{aligned}[t]
           &~\int_{X^\circ} \abs{\nabla^{(0,1)} \zeta}_{\vphilist}^2
           \:\abs{\psi_D}^{1-\eps} +\int_{X^\circ}
           \idxup{\ibddbar\vphi
             +\frac{1-\eps
             }{\abs{\psi_D}} \ibddbar\psi_D }
           \ptinner\zeta\zeta_{\vphilist} \:\abs{\psi_D}^{1-\eps}
           \\
           &
           -2\paren{1-\eps} \Re
           \int_{X^\circ}
           \inner{\dfadj\zeta}{
             \frac{
               \idxup{\diff\psi_D} \ctrt \zeta
             }{
               \abs{\psi_D}
             }
           }_{\mathrlap{\vphilist}} \;\;\:\abs{\psi_D}^{1-\eps}
           +
           \eps
           \int_{X^\circ} \frac{1-\eps }{\abs{\psi_D}^2}
           \abs{\idxup{\diff\psi_D} \ctrt \zeta}_{\vphilist}^2
           \:\abs{\psi_D}^{1-\eps}
         \end{aligned}
    \end{align*}
    for any compactly supported $\zeta \in \smform*/n,q/\paren{X^\circ
      \setminus D; D\otimes F}$, where $\dfadj$ is the formal
    adjoint of $\dbar$ with respect to $\norm\cdot_{\vphi,\omega}$
    (see \cite{Chan&Choi_injectivity-I}*{\S 2.4} for more details).
    Putting $\theta_{\eps'} u$ into $\zeta$ for some $\eps' > 0$
    (assuming that the singularities of $\vphi_F$ has been dealt with
    or $\vphi_F$ is smooth), together with the facts
    \begin{itemize}
    \item \eqref{eq:nable-u=0_ibddbar-vphi_F-u=0} and $\ibddbar\phi_D
      = 0$ on $X \setminus D$,
    \item ($u$ being harmonic) $\dbar\paren{\theta_{\eps'} u} = \dbar\theta_{\eps'} \wedge
      u \;$, $\nabla^{(0,1)}\paren{\theta_{\eps'} u} =
      \dbar\theta_{\eps'} \otimes u$ and
      \begin{equation*}
        \dfadj\paren{\theta_{\eps'} u}
        =-\idxup{\diff\theta_{\eps'}} \ctrt u
        = -\frac{\eps' \theta'_{\eps'}}{\abs{\psi_D}^{1+\eps'}}
        \idxup{\diff\psi_D} \ctrt u \; ,
      \end{equation*}
    \item $\abs{\dbar\theta_{\eps'} \wedge u}_{\vphilist}^2
      +\abs{\idxup{\diff\theta_{\eps'}} \ctrt u}_{\vphilist}^2
      =\abs{\dbar\theta_{\eps'} \otimes u}_{\vphilist}^2$ on $X^\circ$ (see
      \cite{Chan&Choi_injectivity-I}*{footnote 9}),
    \end{itemize}
    one obtains
    \begin{equation*}
      0
      =
      \begin{aligned}[t]
        &-\int_{X^\circ} \frac{1-\eps}{\abs{\psi_D}^\eps} \idxup{\ibddbar\sm\vphi_D}
        \ptinner{\theta_{\eps'} u}{\theta_{\eps'} u}_{\vphilist}
        \\
        &+2(1-\eps) \Re \int_{X^\circ}
        \inner{
          \frac{\eps'\theta'_{\eps'} \idxup{\diff\psi_D} \ctrt u}{\abs{\psi_D}^{1+\eps'}}
        }{
          \frac{\idxup{\diff\psi_D} \ctrt \theta_{\eps'} u}{\abs{\psi_D}}
        }_{\vphilist} \abs{\psi_D}^{1-\eps} \\
        &+\eps \int_{X^\circ} \frac{1-\eps}{\abs{\psi_D}^{1+\eps}}
        \abs{\idxup{\diff\psi_D} \ctrt \theta_{\eps'} u}_{\vphilist}^2
        \; .
      \end{aligned}
    \end{equation*}
    Dividing $(1-\eps)$ from both sides and rearranging terms yield
    \begin{equation*}
      \int_{X^\circ} \frac{\theta_{\eps'}^2}{\abs{\psi_D}^\eps}
      \idxup{\ibddbar\sm\vphi_D}
      \ptinner{u}{u}_{\vphilist}
      =
      \int_{X^\circ} \paren{
        \frac{\eps \theta_{\eps'}^2}{\abs{\psi_D}^{1+\eps}}
        +\frac{2\eps' \theta'_{\eps'} \theta_{\eps'}}{\abs{\psi_D}^{1+\eps+\eps'}}
      }
      \abs{\idxup{\diff\psi_D} \ctrt  u}_{\vphilist}^2 \; .
    \end{equation*}
    Since $\ibddbar\sm\vphi_D$ is a smooth $(1,1)$-form on $X$, the
    integral on the left-hand-side has finite limits as $\eps'
    \tendsto 0^+$ and $\eps \tendsto 0^+$ by the dominated convergence
    theorem.
    As both $\theta_{\eps'}$ and $\theta'_{\eps'}$ are non-negative
    for all $\eps' > 0$, it follows (from letting $\eps' \tendsto
    0^+$) that $\frac{
      \abs{\idxup{\diff\psi_D} \ctrt  u}_{\vphilist}^2
    }{\abs{\psi_D}^{1+\eps}}$ is integrable on $X^\circ$ for any $\eps
    > 0$ and, as a result,
    \begin{equation} \label{eq:limit-of-tBK}
      \int_{X} \idxup{\ibddbar\sm\vphi_D} \ptinner{u}{u}_{\vphilist}
      =\lim_{\eps \tendsto 0^+}
      \eps \int_{\mathrlap{X^\circ}} \;\;
      \frac{
        \abs{\idxup{\diff\psi_D} \ctrt u}_{\vphilist}^2
      }{\abs{\psi_D}^{1+\eps}} 
    \end{equation}
    (see \cite{Chan&Choi_injectivity-I}*{Prop.~3.2.8}; see also
    \cite{Chan&Choi_injectivity-I}*{Lemma 3.3.3} for the treatment
    when $\vphi_F$ is singular).
    
    Via Takegoshi's argument (see
    \cite{Chan&Choi_injectivity-I}*{\S 3.1, Step IV}) 
    and the assumption $[u]\in(\ker\iota_0)^\perp$, one obtains
    \begin{equation}\label{E:Takegoshi_argument}
      \int_{X} \idxup{\ibddbar\sm\vphi_D} \ptinner{u}{u}_{\vphilist}
      =
      \iinner{\ibddbar\sm\vphi_D \Lambda_{\omega}u}{\:u}_{\vphilist} = 0 \; .
    \end{equation}
    On the other hand, since $\frac{u}{\sect_D}$ is smooth on $X^\circ$ and, on any
    admissible open set $V \Subset X$, one has
    \begin{equation*}
      \res{\diff\psi_D}_V = \sum_{j = 1}^{\sigma_V} \frac{dz_j}{z_j}
      -\diff\sm\vphi_D
      \quad\paren{\text{and } \res{\sect_D}_V =z_1 \dotsm z_{\sigma_V}} \; ,
    \end{equation*}
    it then follows from \eqref{E:adjoint_ideal_sheaf_neat_analytic_sing} that
    \begin{equation*}
      \idxup{\diff\psi_D} \ctrt u
      \in \smform/n,{q-1}/[X]\paren{X^\circ ; D \otimes F \otimes
        \aidlof|1|*} \; .
    \end{equation*}
    Proposition \ref{prop:residue-product-X-to-lcS} and
    \eqref{E:Takegoshi_argument} then give
    \begin{equation*}
      0
      =\lim_{\eps \tendsto 0^+} \eps \int_{\mathrlap{X^\circ}\;\;}
      \frac{
        \abs{\idxup{\diff\psi_D}* \ctrt u}_{\vphilist|\phi_D+\vphi_F|}^2
      }{\abs{\psi_D}^{1+\eps}}
      =\pi \sum_{i\in I_D} \int_{D_i}
      \abs{\PRes[D_i](
        \!\idxup{\diff\psi_D}* \ctrt
        \frac{u}{\sect_D}
        )}_{\vphilist F}^2
      \; , 
    \end{equation*}
    which implies that
    \begin{equation*}
      \PRes[D_i]( \idxup{\diff\psi_{D}} \ctrt
      \frac{u}{\sect_{D}} ) \equiv 0 \quad\text{ on } D_i
      \text{ for all } i \in \Iset|| \; .
    \end{equation*}
    Via the argument in \cite{Chan&Choi_injectivity-I}*{\S 3.1, Step V,
        (eq 3.1.5$'$)}, it follows that
    \begin{equation*}
      \text{coef.~of~ }
      \idxup{\diff\psi_{D}} \ctrt \frac{u}{\sect_{D}}
      \text{ ~are locally in }
      \smooth_{X\,*} \text{ on }  X^\circ \; .
    \end{equation*}
    It follows that
    \begin{align*}
      \abs{\:\iinner{\idxup{\diff\psi_D} \ctrt su}{\:
      \frac{\theta'_\eps v_{(\infty)}}{\abs{\psi_D}^{1+\eps}} 
      }_{\vphilist|\varphi| M}} 
      \leq
      &~\int_{\mathrlap{X}\;\;}
        \frac{
        \abs{
        \inner{
        \idxup{\diff\psi_D}* \ctrt
        \frac{u}{\sect_D}
        }{
        \:v_{(\infty)} e^{-\frac 12(\sm\vphi_D+\vphi_M)}}
        }_{\vphilist F}
        \:\abs{s}_{\vphi_M}
        \:\abs{\theta'_\eps}
        }{\abs{\sect_D}_{\sm\vphi_D} \:\abs{\psi_D}^{1+\eps}} \;,
    \end{align*}
    which is finite as $\frac 1{\abs{\sect_D}_{\sm\vphi_D}}$ is
    locally integrable while the other factors are bounded from
    above. 
    This leads to the desired vanishing of the inner product as $\eps
    \tendsto 0^+$, thus completes the proof. \qedhere
  \end{description}
\end{proof}

When $(X,D)$ is plt, i.e.~all $D_i$'s in $D = \sum_{i \in \Iset} D_i$
are the mutually disjoint, one has $\sigma_{\mlc} = 1$ and the
homomorphism $\tau_1$ in \eqref{eq:commut-diagram_sing-Fujino-conj} is
the identity map.
In order to prove the injectivity of $\mu_1$ in this case, it is
enough to show that $\ker\nu_1=\ker\tau_1=0$ by Theorem
\ref{thm:induction-on-Fujino-conj} and Theorem \ref{thm:main-result}.
It follows from \eqref{E:adjunction_residue_sheaf} that the homomorphism $\nu_1$ is reduced to
\begin{multline*}
        \nu_1 \colon
        \bigoplus_{i \in \Iset||} \cohgp q[D_i]{
        K_{D_i} \otimes \res F_{D_i} \otimes
        \mtidlof[D_i]{\res{\vphi_F}_{D_i}}
        } \\
        \longrightarrow~ \bigoplus_{i \in \Iset||} \cohgp q[D_i]{
        K_{D_i} \otimes \parres{F \otimes M}_{D_i} \otimes
        \mtidlof[D_i]{\parres{\vphi_F +\vphi_M}_{D_i}} 
        } \; ,
\end{multline*}
which maps the $i$-th summand to the $i$-th summand via the
multiplication $\otimes \res s_{D_i}$.
The result in \cite{Matsumura_injectivity}*{Thm.~1.3} or the argument given in
\cite{Chan&Choi_injectivity-I}*{\S 3.1} (with $D = 0$) concludes
readily that each $\res{\nu_1}_{D_i}$ is injective.
This recovers the injectivity theorem for plt pairs in
\cite{Matsumura_injectivity-lc}*{Thm.~3.16}.



In our recent preprint with Shin-ichi Matsumura
\cite{Chan&Choi&Matsumura_injectivity}, under the assumption 
that $\varphi_F$ and $\varphi_M$ are smooth, the claims
$\ker\nu_{\sigma} =\ker\tau_{\sigma}$ for $\sigma =1, \dots,
\sigma_{\mlc}$ with general $\sigma_{\mlc} \geq 1$ are proved.
In fact, by considering a commutative diagram analogous to
\eqref{eq:commut-diagram_sing-Fujino-conj}
induced from the short exact sequences
\begin{equation}\label{eq-ex2}
  \xymatrix{
    0 \ar[r]
    & {\faidlof/|\rho|*} \ar[r]
    & {\faidlof|\tau|/|\rho|*} \ar[r]
    & {\faidlof|\tau|/*} \ar[r]
    & 0
  } \quad\text{ for } 0 \leq \rho \leq \sigma \leq \tau \; ,
\end{equation}
together with the isomorphisms $\residlof* \isom \faidlof/-1*$ for
$\sigma = 1, \dots, \sigma_{\mlc}$ obtained from
\eqref{eq:short-ext-seq-of-ideals}, the following
injectivity theorem for a special kind of singular spaces
is obtained.

\begin{thm}[\cite{Chan&Choi&Matsumura_injectivity}*{Thm.~1.2}]
Let $(X,\omega)$ be a compact K\"ahler manifold and $D$ a
reduced snc divisor on $X$.
Let $(F,\varphi_F)$ and $(M,\varphi_M)$ be hermitian line bundles on $X$ equipped with smooth hermitian metrics satisfying 
\begin{equation*}
\ibddbar\vphi_F \geq 0
\;\;\;\text{and}\;\;\;
\ibddbar \vphi_M \leq C \ibddbar \vphi_F \quad\text{ on
} X \;\text{ for some constant } C > 0 \; .
\end{equation*}
Let $s$ be a section of $M$  
such that $s$ does not vanish identically on any lc centres of $(X,D)$.
Then, the multiplication map
\begin{equation*}
H^q(D, K_D \otimes F)
\xrightarrow{\otimes s } 
H^q(D, K_D \otimes F\otimes M)
\end{equation*} 
induced by the tensor product with $s$ is injective for every $q \geq 0$.
\end{thm}

\section{An extension theorem}
\label{sec:extension}

Note that the positivity assumption 
\begin{subequations}
  \begin{equation} \label{eq:pos-assumption-X}
    \ibddbar\paren{\vphi_F +\phi_D +\lambda\psi_D} \geq 0 \quad\text{ for all }
    \lambda \in [0, \delta] \text{ on } X 
  \end{equation}
  for some given constant $\delta > 0$
  is the standard positivity assumption taken in the $L^2$ extension
  theorems (see, for example, \cite{Demailly_extension} and
  \cite{Chan&Choi_ext-with-lcv-codim-1}) and is stronger than the
  assumption in previous sections.
  In what follows, consider also the positivity assumption
  \begin{equation} \label{eq:pos-assumption-lcS}
    \ibddbar\paren{\vphi_F +\phi_{(p)} +\lambda\psi_{(p)}} \geq 0 \quad\text{ for all }
    \lambda \in [0, \delta] \text{ on } \lcS \text{ for each } p \in\Iset
  \end{equation}
  and for each $\sigma =1, 2, \dots, \sigma_{\mlc}$.
\end{subequations}
Write \eqref{eq:pos-assumption-lcS}${}_{\sigma}$ to mean the
assumption for the specific $\sigma$.
It will become apparent that it is convenient to set
\begin{equation*}
  \text{\eqref{eq:pos-assumption-lcS}}_{0} :=
  \text{\eqref{eq:pos-assumption-X}} 
\end{equation*}
and
\begin{equation*}
  \aidlof|-1|* := 0 \quad\text{ and }\quad
  \residlof|0|* :=\faidlof|0|/|-1|* =\aidlof|0|* \; .
\end{equation*}
Consider the long exact sequence
\begin{equation*}
  \xymatrix{
    {\dotsm} \ar[r]
    &{\spH/q-1/{\faidlof|\sigma'|/*}} \ar[r]
    &{\spH{\residlof*}} \ar[r]^-{\tau_\sigma^{\sigma'}}
    &{\spH{\faidlof|\sigma'|/-1*}} \ar[r]
    &{\spH{\faidlof|\sigma'|/*}} \ar[r]
    &{\dotsm}
  }
\end{equation*}
for each $q \geq 1$ and each $\sigma = 0, 1, \dots, \sigma'$.
The following theorem can then be obtained via the techniques
developed in \cite{Chan&Choi_injectivity-I}.

\begin{thm} \label{thm:global-extension-without-estimate}
  Under the positivity assumption
  \eqref{eq:pos-assumption-lcS}${}_{\sigma}$ (for a given $\sigma$),
  one has $\ker\tau_\sigma^{\sigma'} = 0$ for all $q \geq 0$ and for
  any $\sigma' \geq \sigma$.
  The long exact sequence thus splits into short exact sequences 
  \begin{equation*}
    \xymatrix{
      {0} \ar[r]
      &{\spH{\residlof*}} \ar[r]^-{\tau_\sigma^{\sigma'}}
      &{\spH{\faidlof|\sigma'|/-1*}} \ar[r]
      &{\spH{\faidlof|\sigma'|/*}} \ar[r]
      &{0}
    }
  \end{equation*}
  for all $q \geq 0$ and any $\sigma' \geq \sigma$.
  In particular, if the assumption \eqref{eq:pos-assumption-lcS} holds
  for all $\sigma =0, 1, \dots, \sigma'$ for some given $\sigma' \geq
  1$, then, for every 
  \begin{equation*}
  f \in
  \spH/0/{\faidlof|\sigma'|/|\sigma'-1|*} =\cohgp 0[\lcc|\sigma'|']{\logKX
    \otimes \faidlof|\sigma'|/|\sigma'-1|*} \; ,
  \end{equation*}
  there exists a global holomorphic section 
  \begin{equation*}
  F \in \spH/0/{\aidlof|\sigma'|*} =\cohgp
  0[X]{\logKX \otimes \aidlof|\sigma'|*}
  \end{equation*}
  such that $f
  \equiv F \mod \aidlof|\sigma' -1|*$ on $X$, i.e.~$F$ is a
  holomorphic extension of $f$.
\end{thm}

\begin{remark}
  The last claim in Theorem
  \ref{thm:global-extension-without-estimate} is the claim suggested
  in \cite{Chan_Residue-fct-proceedings}*{Conj.~2.2.3}, but with a slightly
  different positivity assumption and without any estimates for the
  extensions.
\end{remark}

\begin{proof}
  Note that $\res{\defidlof{\lcc+1'}}_{\lcS} =\defidlof{\Diff_p D}
  =\genby{\sect_{(p)}} =\mtidlof<\lcS>{\phi_{(p)}}$
  and thus
  \begin{align*}
    K_{\lcS} \otimes \res F_{\lcS} \otimes \mtidlof<\lcS>{\vphi_F}
    &\isom K_{\lcS} \otimes \Diff_p D \otimes \res F_{\lcS} \otimes
    \mtidlof<\lcS>{\vphi_F} \cdot \defidlof{\Diff_p D} \\
    &=K_{\lcS} \otimes \Diff_p D \otimes \res F_{\lcS} \otimes
    \aidlof|0|<\lcS>{\vphi_F}[\psi_{(p)}] \\ 
    &=K_{\lcS} \otimes \Diff_p D \otimes \res F_{\lcS} \otimes
    \mtidlof<\lcS>{\phi_{(p)} +\vphi_F} \; .
  \end{align*}
  It can be seen from this as well as the argument below that the
  computations on $\lcS$ and $X$ are formally the same, so only the
  proof of $\ker\tau_0^{\sigma'} = 0$ is given below.

  Let $\omega$ be a complete K\"ahler form on $X^\circ :=X \setminus
  \vphi_F^{-1}(-\infty)$ as described in
  \cite{Chan&Choi_injectivity-I}*{\S 2.2, item (4)}.
  Let $\abs\cdot_{\phi_D} :=\abs{\cdot}_{\vphi_F+\phi_{D},\omega}$ and
  $\norm\cdot_{X,\phi_D} :=\norm\cdot_{X, \vphi_F+\phi_{D},
    \omega}$ be respectively the pointwise and $L^2$ norms on $X^\circ$
  induced from data shown in the subscript. 
  From the $L^2$ Dolbeault isomorphism (see
  \cite{Matsumura_injectivity}*{Prop.~5.8} and
  \cite{Matsumura_injectivity-lc}*{Prop.~2.8(1)}), one has 
  \begin{align*}
  \cohgp q[X]{K_{X} \otimes F \otimes
    \mtidlof{\vphi_F}} &\isom \cohgp q[X]{\logKX \otimes
    \mtidlof{\phi_{D}+\vphi_F}}
  \\
  &\isom \Harm/q/<X>{D \otimes F},{\vphi_F+\phi_{D}}_{\omega}
  =:\Harm,{\phi_D} \; ,
  \end{align*}
  the space of 
  harmonic forms with respect to $\norm\cdot_{X,\phi_D}$. 
  Using the Bochner--Kodaira formula and the residue computation shown
  in \cite{Chan&Choi_injectivity-I}*{Prop.~3.2.3, Remark 3.2.4,
    Prop.~3.2.5 and Prop.~3.2.8} 
  (as well as \cite{Chan&Choi_injectivity-I}*{Prop.~3.3.2 and Lemma
    3.3.3} to handle the singularities of $\vphi_F$; also
  cf.~\eqref{eq:nable-u=0_ibddbar-vphi_F-u=0} and \eqref{eq:limit-of-tBK}),
  all $u \in \Harm,{\phi_D}$ should satisfy 
  \begin{gather*}
    \nabla^{(0,1)} u = 0 \; , \quad
    \idxup{\ibddbar\vphi_F}\ptinner{u}{u}_{\phi_D} =0 \quad\text{ on }
    X^\circ 
    \\
    \text{and }\quad
    \begin{aligned}[t]
      \pi \sum_{i \in \Iset||} \int_{D_i} \abs{\PRes[D_i](
        \idxup{\diff\psi_{D}} \ctrt \frac{u}{\sect_{D}}
        )
      }_{\vphi_F,\omega}^2
      &=
      \lim_{\eps \tendsto 0^+} \eps
      \int_{\mathrlap{X^\circ}} \;\;
      \frac{\abs{\idxup{\diff\psi_{D}} \ctrt
          u}_{\phi_D}^2}{\abs{\psi_{D}}^{1+\eps}}
      \\
      &=\int_{X^\circ} \idxup{\ibddbar\sm\vphi_{D}} \ptinner u u_{\phi_D} \; .
    \end{aligned}
  \end{gather*}
  However, the positivity assumption \eqref{eq:pos-assumption-lcS}${}_{0}$
  implies that $\idxup{\ibddbar\sm\vphi_{D}} \ptinner u u_{\phi_D} \leq
  0$, which then forces the equality
  \begin{equation*}
    \PRes[D_i]( \idxup{\diff\psi_{D}} \ctrt
    \frac{u}{\sect_{D}} ) \equiv 0 \quad\text{ on } D_i
    \text{ for all } i \in \Iset|| \; .
  \end{equation*}
  Via the argument in \cite{Chan&Choi_injectivity-I}*{\S 3.1, Step V,
    (eq 3.1.5$'$)}, it follows that
  \begin{equation*}
    \text{coef.~of~ }
    \idxup{\diff\psi_{D}} \ctrt \frac{u}{\sect_{D}}
    \text{ ~are locally in }
    \smooth_{X\,*} \text{ on }  X^\circ \; .
  \end{equation*}

  Using the \v Cech--Dolbeault isomorphism, every element in
  $\ker\tau_0^{\sigma'} \subset \Harm,{\phi_D}$ can be written as $\dbar
  v_{(2)} +\dbar v_{(\infty)}$, where $v_{(2)}$ is a $\logKX$-valued
  $(0,q-1)$-form on $X^\circ$ with $L^2$ coefficients with respect to
  $\norm\cdot_{X,\phi_D}$, while $v_{(\infty)}$ is a smooth
  $\logKX$-valued $(0,q-1)$-form on $X^\circ$ which is $L^2$ with
  respect to $\norm\cdot_{X,\vphi_F+\sm\vphi_D}$ (but need not be
  $\norm\cdot_{X,\phi_D}$) (cf.~\cite{Chan&Choi_injectivity-I}*{Lemma
    3.2.1}).
  Let $\theta \colon [0,\infty) \to [0,1]$ be a smooth non-decreasing
  cut-off function and set $\theta_\eps := \theta \circ \frac
  1{\abs{\psi_D}^\eps}$ and $\theta'_\eps := \theta' \circ \frac
  1{\abs{\psi_D}^\eps}$ for $\eps \geq 0$ as in
  \cite{Chan&Choi_injectivity-I}*{\S 3.1} (where $\theta'$ is the
  derivative of $\theta$).
  Recall that $u$ is harmonic.
  Therefore, a computation similar to \cite{Chan&Choi_injectivity-I}*{\S
    3.1, Step III and Step V} yields
  \begin{align*}
    \abs{\iinner{\dbar v_{(2)} +\dbar v_{(\infty)}}{u}_{X,\phi_D}}
    =&~\abs{\iinner{\dbar v_{(\infty)}}{u}_{X,\phi_D}}
    \\
    \xleftarrow{\eps \tendsto 0^+}
     &~\abs{\iinner{\theta_\eps \dbar v_{(\infty)}}{u}_{X,\phi_D}}
    \\
    =
     &~\abs{\cancelto{0}{\iinner{\dbar\paren{\theta_\eps  v_{(\infty)}}}{u}}_{X,\phi_D}
       -\eps \iinner{
       \frac{\theta'_\eps v_{(\infty)}}{\abs{\psi_D}^{1+\eps}}
       }{\: \idxup{\diff\psi_D} \ctrt u }_{X,\phi_D}}
    \\
    \leq
     &~\eps \int_{X^\circ} \frac{\abs{\theta'_\eps} \:\abs{
         \inner{v_{(\infty)} e^{-\frac 12 \sm\vphi_D}} {
           \idxup{\diff\psi_D} \ctrt \frac u{\sect_D}
         }_{\vphi_F,\omega}
       }} {
         \abs{\sect_D}_{\sm\vphi_D} \:
         \abs{\psi_D}^{1+\eps}
       }
    \\
    \xrightarrow{\eps \tendsto 0^+}
    &~0 \; ,
  \end{align*}
  that is, $u \in \paren{\ker\tau_0^{\sigma'}}^\perp$.
  However, $u \in \Harm,{\phi_D}$ is arbitrary, so
  $\ker\tau_0^{\sigma'} = 0$.

  Once the statements $\ker\tau_\sigma^{\sigma'} = 0$ for $\sigma = 0, 1,
  \dots, \sigma'$, and hence the short exact sequences, are obtained,
  the last claim follows from a simple induction.
\end{proof}


\end{document}